\documentclass[journal]{IEEEtran}

\usepackage{graphicx}
\usepackage{color}
\graphicspath{{../pdf/}{../jpeg/}}
\DeclareGraphicsExtensions{.pdf,.jpeg,.png}

\usepackage[cmex10]{amsmath}
\usepackage{amssymb}
\usepackage{array}
\usepackage{mdwmath}
\usepackage{url}
\usepackage{float}
\usepackage{subfig}
\usepackage{multirow}

\begin{document}
\bstctlcite{IEEEexample:BSTcontrol}
    \title{Economic Dispatch Considering Spatial and Temporal Correlations of Multiple Renewable Power Plants}
  \author{Chenghui~Tang\IEEEauthorrefmark{1},~\IEEEmembership{Student Member,~IEEE,}
      Yishen~Wang\IEEEauthorrefmark{2},~\IEEEmembership{Student Member,~IEEE,} \\
      Jian Xu\IEEEauthorrefmark{1},~\IEEEmembership{Member,~IEEE,}
      Yuanzhang~Sun\IEEEauthorrefmark{1},~\IEEEmembership{Senior Member,~IEEE}
      and~Baosen~Zhang\IEEEauthorrefmark{2},~\IEEEmembership{Member,~IEEE,}\\

\thanks{This work was supported in part by the National Key R\&D Program of China (2016YFB0900105), in part by the National Natural Science Foundation of China (51477122)}
\thanks{
	C. Tang, J. Xu and Y. Sun are with the School of Electrical Engineering, Wuhan University, Wuhan, 430072 China.}
\thanks{
	Y. Wang and B. Zhang are with the Department of Electrical Engineering, University of Washington, Seattle, WA, 98195 USA.}}

\maketitle


%

\begin{abstract}
	The correlations of multiple renewable power plants (RPPs) should be fully considered in the power system with very high penetration renewable power integration. This paper models the uncertainties, spatial correlation of multiple RPPs based on Copula theory and actual probability historical histograms by one-dimension distributions for economic dispatch (ED) problem. An efficient dynamic renewable power scenario generation method based on Gibbs sampling is proposed to generate renewable power scenarios considering the uncertainties, spatial correlation and variability (temporal correlation) of multiple RPPs, in which the sampling space complexity do not increase with the number of RPPs. Distribution-based and scenario-based methods are proposed and compared to solve the real-time ED problem with multiple RPPs. Results show that the proposed dynamic scenario generation method is much more consist with the actual renewable power. The proposed ED methods show better understanding for the uncertainties, spatial and temporal correlations of renewable power and more economical compared with the traditional ones.
\end{abstract}


%
\IEEEpeerreviewmaketitle

\section{Introduction}
The increase in penetration of renewable resources around the world is fundamentally changing how power systems are operated \cite{Risk_Based_UC}. The potential risks of power unbalance, transmission congestion become much more complicated due to the dependence of wind power plants (WPPs) and solar photovoltaic power plants (PVPPs) in the system. In power system economic dispatch (ED), system operators should fully consider the uncertainties, spatial and temporal correlations of renewable power plants (RPPs). In this paper, we focus on the problem of economic dispatch (ED) when there are multiple correlated RPPs in the system.

The first step in integrating RPPs into dispatch is to obtain a tractable model that captures the uncertainties and correlations between them. Conditional distribution models for RPPs have proved to be a reliable mathematical method for ED with renewable power integration \cite{Versatile}\cite{sce_generation_Ma}\cite{copula_Zhang}. Conditional distribution models employ the forecast power of RPPs to increase the representation accuracy of the uncertainties and correlations of renewable power and consist with the current ED mode. Forecast bins are employed in \cite{Versatile} and \cite{sce_generation_Ma} and then actual wind power in each forecast bin can be modeled by mathematical distribution \cite{Versatile} or actual probability density histogram (PDH) \cite{sce_generation_Ma}. However, when the number of RPPs increases, forecast bin method becomes hard to be employed due to the curse of dimensionality caused by the number of bins. Copula theory shows better potential in dealing with multiple renewable power variables as in \cite{copula_Zhang}. Wind power conditional distribution model is built in \cite{copula_Zhang} based on copula theory and wind power scenarios are generated by the conditional joint distribution to solve the unit commitment and ED. However, high-dimension distribution model would greatly increase the computation scale and even hard to be employed directly. In addition, although wind power uncertainties and spatial correlation are novelly considered in each time interval, the temporal correlation (variability) of wind power scenarios among the schedule horizon which is regarded same important  \cite{sce_generation_Ma}\cite{PCA} is not considered in \cite{copula_Zhang}. In this paper, we do not use the high-dimension conditional joint distribution model directly. Instead, we convert the scenarios generation using the conditional joint distribution model at one time interval (static scenarios) to scenarios generation using the one-dimension conditional distribution by Gibbs theory \cite{Gibbs}. Then, a dynamic scenario generation method is proposed and the renewable power variability is considered among the schedule horizon. In addition, in order to study the overall effect of renewable power uncertainties and correlations on the system risk of reserve deficiency and transmission congestion, we also build the conditional marginal distribution of sum actual power of all RPPs and actual renewable power in power system transmission lines.

The next step is to incorporate the stochastic model into ED in a computationally efficient fashion. Stochastic ED \cite{chance_constrain}\cite{ST_ED} and robust ED \cite{Robust_1}\cite{Robust_2} are two main methods to solve the ED problem with renewable power integration. Compared with robust ED, stochastic ED offers better mechanisms to manage the uncertainties explicitly \cite{ST_ED}. One way to account for uncertainties in stochastic ED is to employ chance constraints to maintain a predefined risk level for the whole system \cite{chance_constrain}\cite{ST_ED}. However, in a system with possibly congested transmission lines and multiple RPPs, it becomes difficult to define a single risk level for the entire system. Here we consider risk by explicitly considering load shedding (LS) and renewable energy curtailment (REC) caused by system reserve deficiency and transmission congestions. To deal with the uncertainties, correlations and variability more reasonably, we seek for the optimal level of risk. To represent the uncertainties, spatial and temporal correlations of RPPs, distribution-based ED method and scenario-based ED method are proposed and compared to solve the real-time economic dispatch (RTED) problem with multiple RPPs. To test the ED results, we use the scenarios by proposed dynamic scenario generation method. The main contributions of this paper are summarized as follows.



We propose an efficient dynamic scenario generation method that captures the joint distribution and the variability
of multiple renewable power plants. Based on Gibbs sampling, our proposed method avoids directly computing high dimensional distributions to greatly reduce the complexity of sampling correlated renewable power plants.

We then describe two ED methods: distribution-based and scenario-based and compare their performances on economic dispatch problem with multiple renewable power plants. The potential risk of load shedding and renewable power curtailment caused by uncertainties and correlations of renewable power integration are modeled to balance the conventional power plants (CPPs) outputs, system reserve, potential risk of load shedding and renewable power curtailment.

We find that the best method to use depends on the number of scenarios used. With a small set of scenarios, scenario-based economic dispatch shows much worse compared with of the distribution-based one. With the increase of the number of scenarios embedded in the economic dispatch, scenario-based economic dispatch performs better. We provide detailed discussion of the results and illustrate which one should be chosen in practice based on the computational power and information available to the system operator. 

The remaining of this paper is organized as follows. In Section II, uncertainties and spatial correlation of multiple RPPs are modeled by Copula theory and actual PDH. In Section III, static scenario generation method of multiple RPPs based on Gibbs theory is proposed. Dynamic renewable power scenarios are generated to represent the uncertainties, spatial and temporal correlations of RPPs. A distribution-based and scenario-based RTED method are proposed in Section IV and Section V, respectively. In Section VI, uncertainties and correlations of multiple RPPs are shown and the dynamic renewable power scenarios are discussed. Numerical experiments of the proposed RTED models are conducted and compared using the IEEE 118-bus system. Section VII provides conclusions.

\section {Modeling the Uncertainties and Spatial Correlation of Multiple RPPs in Power System}
The uncertainties and spatial correlation of renewable power plants (RPPs) should be modeled to evaluate the potential risk of the load shedding {(LS)} and renewable energy curtailment {(REC)}. This section describes how uncertainties and spatial correlation of multiple RPPs are modeled based on historical data (forecast and actual power of each RPP) using copula theory.
By uncertainty, we mean the marginal distribution of the forecast error of each RPP; and by spatial correlation, we mean the joint relationship between the forecast errors. Note that we often interested in the \emph{conditional distribution} of the errors given the forecasted values.


\subsection {Conditional Distribution of Renewable Power}
The power productions of the RPPs are described by two random vectors: a vector of \emph{forecast values} and a vector of \emph{actual power production} conditioned on the forecasts. Let $w_{f,j}$ denote the forecast of the $j$'th RPP if it is a wind power plant (WPP) and $s_{f,k}$ denote its forecast if it is solar photovoltaic power plant (PVPP). We assume there are $J$ total WPPs and $K$ total PVPPs. Let {$\mathbf{f}$} denote the vector of forecasts, i.e., $(w_{f,{1}}...w_{f,J},s_{f,{1}}...s_{f,K})$. Let $F(w_{f,i})$ denote the marginal cumulative distribution function (CDF) of the forecast of the $i$'th WPP (similar for PVPPs). We use {$\Omega(\mathbf{f})$} to denote the set of marginal forecast CDFs, i.e., $(F(w_{f,{1}})...F(w_{f,J}),F(s_{f,{1}})...F(w_{s,K}))$.

Copula method is an effective way of modeling the multiple variables dependence where the marginal distributions and the correlations between random variables appear separately~\cite{copula_Zhang,WytockEtAl2013}.
In the Copula method, the joint CDF of the forecasted and actual RPP productions, $F(w_{a,{1} }. . .w_{a,J},s_{a,{1}}. . .s_{a,K},\mathbf{f})$, is written as:
\begin{equation} \label{eq:01}
\begin{aligned}
& F(w_{a,{1} }. . .w_{a,J},s_{a,{1}}. . .s_{a,K},\mathbf{f}) \\
= & C(F(w_{a,{1}})...F(w_{a,J}),F(s_{a,{1}})...F(s_{a,K}), \Omega(\mathbf{f}))
\end{aligned}
\end{equation}
where the function $C(\cdot)$ is called the Copula function. Essentially, we transform the joint CDF  $F(w_{a,{1} }. . .w_{a,J},s_{a,{1}}. . .s_{a,K},\mathbf{f})$ to a function of the marginal CDFs  $F(w_{a,j})$, $F(s_{a,k})$, $F(w_{f,j})$, $F(s_{f,k})$ linked by the Copula function $C$. Similarly, the joint probability density function (PDF) is:
\begin{equation} \label{eq:02}
\begin{aligned}
& f(w_{a,{1} }. . .w_{a,J},s_{a,{1}}. . .s_{a,K},\mathbf{f}) \\
= & c(F(w_{a,{1}})...F(w_{a,J}),F(s_{a,{1}})...F(s_{a,K}), \Omega(\mathbf{f}))\\
& \cdot \prod_{j=1}^{J}f(w_{a,{j}}) \cdot \prod_{k=1}^{K}f(s_{a,{k}}) \cdot \prod_{j=1}^{J}f(w_{f,{j}}) \cdot \prod_{k=1}^{K} f(s_{f,{k}}).
\end{aligned}
\end{equation}
If only the forecasted values are considered, then similar to \eqref{eq:02}, their joint PDF can be written as:
\begin{equation} \label{eq:03}
f(\mathbf{f})= c(\Omega(\mathbf{f})) \cdot \prod_{j=1}^{J}f(w_{f,{j}}) \cdot \prod_{k=1}^{K} f(s_{f,{k}}).
\end{equation}
Lastly, since the forecasted values are already know in ED, all the randomness are left in the actual power of the RPPs given the forecasted values. Combining \eqref{eq:02} \eqref{eq:03}, the joint conditional PDF of the actual productions is:
\begin{equation} \label{cjdistribution}
\begin{aligned}
& f(w_{a,{1} }. . .w_{a,J},s_{a,{1}}. . .s_{a,K}|\mathbf{f}) \\
=& \frac{f(w_{a,{1} }. . .w_{a,J},s_{a,{1}}. . .s_{a,K},\mathbf{f})}{f(\mathbf{f})}\\
=& \frac{c(F(w_{a,{1}})...F(w_{a,J}),F(s_{a,{1}})...F(s_{a,K}),\Omega(\mathbf{f}))}{c(\Omega(\mathbf{f}))} \\
& \cdot \prod_{j=1}^{J}f(w_{a,{j}}) \cdot \prod_{k=1}^{K}f(s_{a,{k}})
\end{aligned}
\end{equation}

 There are many suitable copula functions (e.g., Gaussian, t, empirical~\cite{G_copula}) that can be used in \eqref{cjdistribution}. In this paper, we adopt the Gaussian copula and use the actual {probability density histogram} (PDH) \cite{sce_generation_Ma} to get the marginal distribution functions in Copula method.

\subsection {Sampling the Conditional Distribution}
In Gibbs sampling theory, conditional distribution function of actual available power of each RPP is needed and can be modeled in \eqref{ccdistribution} (take $j$'th WPP for instance), shown in the top of the next page. Note that the conditioning is on the forecast and actual renewable powers except for $j$'th WPP. This model will be further discussed in Section III.
\newcounter{TempEqCnt}
\setcounter{TempEqCnt}{\value{equation}}
\setcounter{equation}{4}
\begin{figure*}[t]
	\hrulefill
	\begin{equation} \label{ccdistribution}
	\begin{aligned}
	f(w_{a,{\it j} }|w_{a,{1} }. . .w_{a,j-1},w_{a,j+1}...w_{a,J},s_{a,{1}}. . .s_{a,K},\mathbf{f}) \\
	=\frac{c(F(w_{a,{1}})...F(w_{a,J}),F(s_{a,{\it 1}})...F(s_{a,K}),\Omega(\mathbf{f}))}{c(F(w_{a,{1}})...F(w_{a,j-1}),F(w_{a,j+1})...F(w_{a,J}),F(s_{a,{1}})...F(s_{a,K}),\Omega(\mathbf{f}))} & \cdot f(w_{a,{\it j}})
	\end{aligned}
	\end{equation}
\end{figure*}
\setcounter{equation}{\value{TempEqCnt}}


\subsection {Sum Power of RPPs}

If congestion is neglected in the ED model, then we only need to model the sum of actual powers of RPPs, denoted by $R_a^{\Sigma}$.
Similar to \eqref{eq:01}, the joint CDF can be written using a Copula function as:
\setcounter{equation}{5}
\begin{equation} \label{eq:3}
\begin{split}
& F(R_{a}^{\Sigma},\mathbf{f})=C(F(R_{a}^{\Sigma}),\Omega(\mathbf{f}))\\
& R_{a}^{\Sigma}=\sum _{{\it j}={1}}^{{\it J}}w_{a,j}+\sum _{{\it k}={1}}^{{\it K}}s_{a,k} \\
\end{split}
\end{equation}
The conditional distribution of sum actual available power of all RPPs given the forecasts is:
\begin{equation} \label{sumdistribution}
f(R_{a}^{\Sigma}|\mathbf{f})= \frac{c(F(R_{a}^{\Sigma}),\Omega(\mathbf{f}))}{c(\Omega(\mathbf{f}))}\cdot   f(R_{a}^{\Sigma})
\end{equation}

\subsection {Conditional Distribution Under Congestion}
The probability of congestions can be evaluated via {shift factors} by finding the contribution of line flows from each RPP~\cite{li2013dynamic}.

Let $R_{a}^{L_l}$ denote the renewable power flow in the transmission line $L_l$.
Based on the synchronous historical forecast and actual power of each RPP, the renewable power in the transmission lines can again be modeled with a copula:
\begin{equation} \label{eq:5}
\begin{aligned}
& F(R_{a}^{L_l},\mathbf{f})= C(F(R_{a}^{L_l}),\Omega(\mathbf{f}))\\
& R_{a}^{L_l}=\sum_{j=1}^{J}k_{l,j}w_{a,j}+\sum_{k=1}^{K}k_{l,k}s_{a,k} \\
\end{aligned}
\end{equation}
{where $k_{l,j}$ and $k_{l,k}$ are the generation distribution shift factors of $j$'th WPP and $k$'th PVPP to transmission line $L_l$, respectively.}
Based on the forecast power of RPPs $w_{f,j}$ and $s_{f,k}$, the conditional distribution of renewable power in the transmission lines $RE_{a}^{L_l}$ can be modeled as follows.
\begin{equation} \label{linedistribution}
\begin{aligned}
f(R_{a}^{L_l}|\mathbf{f})
= \frac{c(F(R_{a}^{L_l}),\Omega(\mathbf{f}))}{c(\Omega(\mathbf{f}))} \cdot f(R_{a}^{L_l})
\end{aligned}
\end{equation}

In this paper, (5), (7) and (9) are used in the RTED model. By the suitable Copula function and actual PDH, we can model the renewable power uncertainties and spatial correlation accurately using just one-dimension distributions.



\section{Scenario Generation}
In this section, we first propose a reliable static renewable power scenario generation method in each time interval $1,\dots,T$. Then we present an efficient dynamic renewable power scenario generation method for the entire time horizon.

\subsection {Static Scenario Generation}

By the joint distribution of multiple RPPs in \eqref{cjdistribution}, scenarios can be generated to represent the uncertainties and spatial correlation of all RPPs in the system. However, with the increase of the number of RPPs, classical random sampling methods such as inverse transform sampling and Latin hypercube sampling \cite{L_sampling} become hard to be employed due to matrix size and computational limitations. Other classical sampling methods such as rejection sampling tend to have very large rejection rate for a high number of dimensions.

To this end, a reliable static renewable power scenario generation method based on Gibbs sampling \cite{Gibbs} is proposed to sample for the conditional joint distribution function of actual available power of RPPs in \eqref{cjdistribution}. Compared with directly sampling by the conditional joint distribution \cite{copula_Zhang}, Gibbs sampling converts the sampling process of joint distribution in \eqref{cjdistribution} to $J+K$ sampling processes of conditional distribution in \eqref{ccdistribution}. Namely, let $U$ be a random variable generated uniformly within $[0,1]$, then each RPP can be sampled via the inverse transform:
\begin{equation} \label{inversesampling}
w_{a,j}=F_{a,j}^{-1}(U),\quad s_{a,k}=F_{a,k}^{-1}(U)
\end{equation}
where $F_{a,j}^{-1}$ and $F_{a,k}^{-1}$ is the inverse function of $F_{a,j}$ and $F_{a,k}$, respectively.

Gibbs sampling needs a burn-in process \cite{burn_in} before it converges to the true distribution in \eqref{cjdistribution}. So we throw out $N_{b}$ (e.g. 1000) samples in the beginning the process. The detailed procedure of static scenarios generation is:
\begin{enumerate}
	\item Setting the number of renewable power scenarios: $N_{sc}$ (e.g. 5000), the total number of samples is $N_{sc}+N_{b}$.
	\item Setting the initial sampling values to be the forecasted power for each RPP.
	\item Employing inverse transform sampling in \eqref{inversesampling} in a round robin fashion for each scenario generation step (indexed by $i$):

\begin{itemize}
	\item $f(w_{a,{1}}^{i}|w_{a,2}^{i}...w_{a,J}^{i},s_{a,{1}}^{i}...s_{a,K}^{i},\mathbf{f})$
	\item $f(w_{a,{\it j}}^{i}|w_{a,{1}}^{i+1}...w_{a,{{\it j}-1}}^{i+1},w_{a,{{\it j}+1}}^{i}...w_{a,J}^{i},s_{a,{1}}^{i}...s_{a,K}^{i},\mathbf{f})$
	\item $...$
	\item $f(s_{a,{\it k}}^{i}|w_{a,{1}}^{i+1}...w_{a,J}^{i+1},s_{a,{1}}^{i+1}...s_{a,{{\it k}-1}}^{i+1},s_{a,{{\it k}+1}}^{i}...s_{a,K}^{i},\mathbf{f})$
	\item $f(s_{a,{\it K}}^{i}|w_{a,{1}}^{i+1}...w_{a,J}^{i+1},s_{a,{1}}^{i+1}...s_{a,{{\it K}-1}}^{i+1},\mathbf{f})$
\end{itemize}

	\item Repeating 3 from {\it i}=1...$N_{sc}+N_{b}$. Disregard the first $N_{b}$ scenarios and we get $N_{sc}$ renewable power scenarios.

\end{enumerate}

{An important feature of the proposed static scenario generation method is that with the increase of the number of RPPs, the computational space complexity remains same and the computational time complexity increases linearly, effectively mitigating the curse of dimensionality.}

\subsection {Dynamic Scenario Generation}
{A dynamic scenario is a scenario that considers the variability (i.e., temporal correlation) of the output of a RPP.} The method presented in the last section can generate renewable power scenarios of conditional joint distribution (c.f. \eqref{cjdistribution}) which captures the marginal uncertainties and spatial correlation. In this section we extend it to capture the temporal correlation among the time points in a scenario, which is also of vital importance in power system operations~\cite{sce_generation_Ma,PCA,sce_generation_Pinson}.

To capture the variability, some new variables are introduced. Take a WPP for instance, a new random variable $Z_{a,j}^{t}$ is introduced which follows
the standard Gaussian distribution with zero mean and unit standard deviation. Since the value of CDF of $Z_{a,j}^{t}$ is uniformly distributed over [0,1], the uniform distribution $U$ in \eqref{inversesampling} can be replaced by a CDF $\Phi(Z_{a,j}^{t})$.  Given the realization of random variable $Z_{a,j}^{t}$, $w_{a,j}^{t}$ can be sampled as follows:

\begin{equation} \label{transform}
\begin{aligned}
w_{a,j}^t=F_{a,j}^{-1}(\Phi(Z_{a,j}^{t}))
\end{aligned}
\end{equation}

To consider the variability of each RPP, it is assumed that the joint distribution of $Z_{a,j}^{t}$ follows a multivariate Gaussian distribution $Z_{a,j}^{t} \sim N(\mu_{j},\Sigma_{j})$. The expectation of $\mu_{j}$ is a vector of zeros and the covariance matrix $\Sigma_{j}$ satisfies

\begin{equation} \label{matrix}
\Sigma_j=\left[
\begin{matrix}
\sigma_{1,1}^{j}&\sigma_{1,2}^{j}&\dots&\sigma_{1,{\it T}}^{j}&\\
\sigma_{2,1}^{j}&\sigma_{2,2}^{j}&\dots&\sigma_{2,{\it T}}^{j}&\\
\vdots&\vdots&\ddots&\vdots&\\
\sigma_{{\it T},1}^{j}&\sigma_{{\it T},2}^{j}&\dots&\sigma_{{\it T},{\it T}}^{j}&\\
\end{matrix}
\right]
\end{equation}

\noindent where $\sigma_{m,n}^{j}=cov(Z_{a,j}^{m},Z_{a,j}^{n})$, {\it m}, {\it n}=1,2...{\it T}, $\sigma_{{\it m}, {\it n}}^{j}$ is the covariance of $Z_{a,j}^{m}$ and $Z_{a,j}^{n}$.

\begin{figure}[!htb]
	\begin{center}
		\includegraphics[trim = 10 250 60 200, clip, width=1.0\columnwidth]{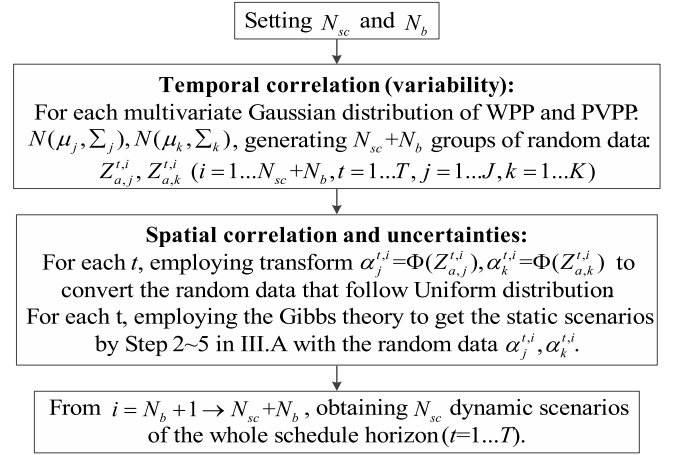}\\
		\caption{Flowchart of dynamic renewable power scenario generation method}\label{flowchart}
	\end{center}
\end{figure}

The covariance structure of $\Sigma_j$ can be identified by covariance $\sigma_{m,n}^{j}$. As is done in \cite{sce_generation_Ma}\cite{sce_generation_Pinson}, an exponential covariance function is employed to model $\sigma_{m,n}^{j}$ in \eqref{matrix},

\begin{equation} \label{exponential}
\begin{aligned}
\sigma_{m,n}^{j}=\rm exp(-\frac{|{\it m}-{\it n}|}{\epsilon_{\it j}}) \quad 0 \le {\it m},  {\it n} \le {\it T}
\end{aligned}
\end{equation}

\noindent where $\epsilon_{\it j}$ is the range parameter controlling the strength of the
correlation of random variables $Z_{a,j}^{t}$ among the set of lead-time. Similar to \cite{sce_generation_Ma}, $\epsilon_{\it j}$ can be determined by comparing the distribution of renewable power variability of the generated scenarios by the indicator in \cite{sce_generation_Ma}. Here, assuming that the  range parameter $\epsilon_{\it j}$ of each RPP have been obtained, the flowchart of dynamic renewable power scenario generation method is as shown in Fig.~\ref{flowchart}.

Before generating $N_{sc}$ scenarios, small amount of scenarios are generated to obtain the range parameter of each RPP. After all the range parameters in \eqref{matrix} are obtained, we can start the dynamic wind power scenarios generation in Fig.~\ref{flowchart}. At each time interval, they follow the conditional joint distribution in \eqref{cjdistribution} and among the time horizon, the variability is considered.

One thing that need to be noticed is that each static scenario generation process in Fig. 1 does not affect each other after the random data set is determined. Parallel computing can be employed to increase the computation efficiency to meet the real-time requirement.

In scenario-based method, the above generated scenarios should be reduced to certain number of scenarios that deemed as the most probability occur. A scenario reduction method in \cite{YishenWang} is employed in this paper for the reason that it has great efficiency compared with other methods to meet the real-time requirement.

\section{Distribution-Based ED} \label{sec:dist}
In this section we study the RTED problem where CPPs outputs, system reserve and potential risk of {renewable energy curtailment} and {load shedding} are balanced. We consider an hourly dispatch with $T=12$ intervals where each one is 5 minutes long. The objective function of the ED problem is:
%
%
\begin{equation} \label{DB_1}
\begin{aligned}
\begin{aligned}
min\sum _{{\it t}={\it 1}}^{{\it T}}E[f_t]& =\sum _{{\it t}=1}^{{\it T}}E[f_{c,t}(p_{i,t},r_{u,i,t},r_{d,i,t})]\\
& +\sum _{{\it t}=1}^{{\it T}}E[f_{R,t}(w_{c,j,t},s_{c,k,t},l_{s,b,t})]\\
\end{aligned}
\end{aligned}
\end{equation}
where $f_t$  is the total system cost at time {\it t};
$f_{c,t}$ is the total CPP cost at time {\it t}; $f_{R,t}$ is the total penalty cost caused by renewable power uncertainties (REC and LS); $p_{i,t}$ is the schedule power of $i$'th CPP at time {\it t}; $r_{u,i,t}$ and $r_{d,i,t}$ is the upward and downward reserve of $i$'th CPP at time {\it t}, respectively; $w_{c,j,t}$ and $s_{c,k,t}$ is the power of REC of $j$'th WPP and $k$'th PVPP at time {\it t}, respectively; $l_{s,b,t}$ is the power of LS of $b$'th bus at time {\it t}.

The CPP cost is given by
\begin{equation} \label{DB_2}
\begin{aligned}
& f_{c,t}(p_{i,t},r_{u,i,t},r_{d,i,t})\\
=\sum _{{\it i}=1}^{{\it I}}(b_{f,i}p_{i,t} & +c_{f,i}+c_{ur,i}r_{u,i,t}+c_{dr,i}r_{d,i,t})
\end{aligned}
\end{equation}
where {\it I} is the total number of CPPs; $b_{f,i}$ and $c_{f,i}$ are the fuel cost coefficients of $i$'th CPP, respectively; $c_{ur,i}$ and $c_{dr,i}$ are the cost coefficients of upward and downward reserve of $i$'th CPP, respectively.
Penalties with respect to uncertainties in the renewable powers are given by:
\begin{equation} \label{DB_3}
E[f_{R,t}(w_{c,j,t},s_{c,k,t},l_{s,b,t})]=c_{ls}E_{ls,t}+c_{rec}E_{rec,t}\\
\end{equation}
where $c_{ls}$ and $c_{rec}$ is the penalty coefficients of LS and REC, respectively; $E_{ls,t}$ and $E_{rec,t}$ is the expected values of LS and REC, respectively.

For ease of analysis, the sum scheduled renewable energy $R_{t}^{\Sigma}$ is introduced as an internal variable in the distribution-based ED model for the balance of power system.


$R_{a,t}^{\Sigma}$ is the sum actual available power of all RPPs at time {\it t} as shown in \eqref{sumdistribution}. $\underline{R}_{t}$ and $\overline{R}_{t}$ is the lower and upper bound that renewable power can be compensated by system reserves at time {\it t}, respectively. In worse case, if the sum actual renewable power locates in the outside of $[\underline{R}_{t}, \overline{R}_{t}]$, system reserve cannot cover all the uncertainties of renewable power. At this time, LS or REC would be employed for the power balance of the system. Then the total penalty cost of renewable power $f_{R,t}(w_{c,j,t},s_{c,k,t},l_{s,b,t})$ can be converted to $f_{R,t}(\underline{R}_{t},\overline{R}_{t})$ and written as
\begin{equation} \label{DB_4}
\begin{aligned}
\begin{aligned}
& E[f_{R,t}(w_{c,j,t},s_{c,k,t},l_{s,b,t})]=f_{R,t}(\underline{R}_{t},\overline{R}_{t})\\
& =c_{ls}\int_{0}^{\underline{R}_{t}}(\underline{R}_{t}-R_{a,t}^{\Sigma})f(R_{a,t}^{\Sigma})dR_{a,t}^{\Sigma}\\
& +c_{rec}\int_{\overline{R}_{t}}^{R_{r}}(R_{a,t}^{\Sigma}-\overline{R}_{t})f(R_{a,t}^{\Sigma})dR_{a,t}^{\Sigma}\\
\end{aligned}
\end{aligned}
\end{equation}
where $R_{r}$ is the total capacity of renewable power.

Compared with other classical stochastic ED methods that use a predefined confidence level to convert the reserve chance constraints to be linear ones \cite{Versatile}\cite{chance_constrain}\cite{ST_ED}, we can seek the optimal confidence level to find the balance for CPPs outputs, system reserve and potential risk of REC and LS according to different situations.


All constraints of the proposed distribution-based ED model are as follows:
\begin{equation} \label{DB_5}
\sum_{{\it i}=1}^{{\it I}}p_{i,t}+R_{t}^{\Sigma}=L_{t} \quad \forall t
\end{equation}
\begin{equation} \label{DB_6}
\begin{aligned}
& R_{t}^{\Sigma}-\sum_{{\it i}=1}^{{\it I}}r_{u,i,t}=\underline{R}_{t} \quad \forall t \\
& R_{t}^{\Sigma}+\sum_{{\it i}=1}^{{\it I}}r_{d,i,t}=\overline{R}_{t} \quad \forall t
\end{aligned}
\end{equation}
\begin{equation} \label{DB_7}
\begin{aligned}
p_{i,t} + r_{u,i,t} \leq p_{max,i} \quad \forall i,t \\
p_{i,t} - r_{d,i,t} \geq p_{min,i} \quad \forall i,t
\end{aligned}
\end{equation}
\begin{equation} \label{DB_8}
\begin{aligned}
p_{i,t}-p_{i,t-1} & \leq \Delta p_{u,max,i} \quad \forall i,t \\
p_{i,t-1}-p_{i,t} & \leq \Delta p_{d,max,i} \quad \forall i,t
\end{aligned}
\end{equation}
\begin{equation} \label{DB_9}
\begin{aligned}
0 & \leq r_{u,i,t} \leq r_{u,max,i} \quad \forall i,t \\
0 & \leq r_{d,i,t} \leq r_{d,max,i} \quad \forall i,t
\end{aligned}
\end{equation}
\begin{equation} \label{DB_10}
\begin{aligned}
0 \leq \underline{R}_{t},\;  \overline{R}_{t}\leq R_{r} \quad \forall t
\end{aligned}
\end{equation}


\vspace{-1em}

\begin{equation} \label{DB_12}
\begin{aligned}
\begin{aligned}
& \sum_{i=1}^{I}k_{l,i}p_{i,t}+\underline{R}_{a}^{L_{\it l}}-\sum_{{\it b}=1}^{{\it Nb}}k_{l,b}L_{b,t} \ge -Pl_l^{max}  \quad \forall l,t\\
& \sum_{i=1}^{I}k_{l,i}p_{i,t}+\overline{R}_{a}^{L_{\it l}}-\sum_{{\it b}=1}^{{\it Nb}}k_{l,b}L_{b,t} \le  Pl_l^{max}  \quad \forall l,t
\end{aligned}
\end{aligned}
\end{equation}
where
\begin{itemize}

	\item \eqref{DB_5} is the supply-demand balance constraint; $L_{t}$ is the forecast power demand at time {\it t};

	\item \eqref{DB_6} is the system reserve constraint;

	\item \eqref{DB_7} are the CPPs scheduled power plus reserve capacity constraint;  $p_{max,i}$  and   $p_{min,i}$ are the upper and lower generation limit of the $i$'th CPP, respectively;

	\item \eqref{DB_8} are the CPPs ramp-rate constraint; $\Delta p_{u,max,i}$ and $\Delta p_{d,max,i}$ are the maximum amount of upward and downward ramp rate of $i$'th CPP within a specific time period (e.g., 5min), respectively;

	\item \eqref{DB_9} are the reserve capacity constraints; $r_{u,max,i}$ and $r_{d,max,i}$ are the maximum amount of up and down reserves that the $i$'th CPP is capable of providing, respectively;

	\item \eqref{DB_10} are the confidence level bound constraint;


	\item \eqref{DB_12} are the transmission capacity constraint; {\it Nb} is the total number of buses; $Pl_l^{max}$ is the transmission capacity limit on transmission line $l$; based on the the distribution of renewable power in the transmission lines $R_{a}^{L_l}$ in \eqref{linedistribution}, the uncertainties and correlations of multiple power energy can be considered compared with the classical model in \cite{AI} and  \cite{VersatileMixture} that used the forecast or scheduled renewable power. A conservative bound such as 99.9\% can be used in this constraint.
\end{itemize}


\section{Scenario-Based ED}
Different from the distribution-based ED, scenario-based ED incorporate the renewable power uncertainties by a certain number of possible renewable power series (i.e., scenarios). This means that scenario-based ED is essentially a deterministic optimization. This allows a more flexible way to model the risk of renewable power such as REC caused by certain transmission line congestion. However, the performance of scenario-based ED greatly relies on the number of scenarios that are considered in the ED. RTED model based on the scenario of multiple RPPs is proposed in this section. The potential risk of LS and REC caused by system reserve deficiency and transmission congestion are modeled by the scenario-based ED. The penalty cost caused by renewable power uncertainties (REC and LS) in (14) and (16) can be written using scenarios as:
\begin{equation} \label{SB_1}
\begin{aligned}
\begin{aligned}
&E[f_{R,t}(w_{c,j,t},s_{c,k,t},L_{s,b,t})]  \\
&=\sum_{{\it sc}=1}^{{\it SC}}[p^{sc}(c_{rec}(\sum_{{\it j}=1}^{{\it J}}w_{c,j,t}^{sc}+\sum_{{\it k}=1}^{{\it K}}s_{c,k,t}^{sc})
+c_{ls} \sum_{{\it b}=1}^{{\it Nb}}L_{s,b,t}^{sc})]
\end{aligned}
\end{aligned}
\end{equation}
where $sc$ is the {\it sc}-th scenario for WPPs and PVPPs, $SC$ is the number of renewable power scenarios in RTED model, $w_{c,j,t}^{sc}$ is the amount of wind power curtailment of $j$'th WPP at time {\it t} of $sc$'th scenario; $s_{c,k,t}^{sc}$ is the amount of solar power curtailment of $k$'th PVPP at time {\it t} of $sc$'th scenario; $L_{s,b,t}^{sc}$ is the amount of LS of $b$'th bus at time {\it t} of $sc$'th scenario.

Then the optimization problem is same as in Section~\ref{sec:dist}, except the constraints and objectives are represented with scenarios. In particular, the constraints are:
\begin{equation} \label{SB_2}
\begin{aligned}
0  \leq w_{c,j,t}^{sc} \leq w_{a,j,t}^{sc} \quad \forall j,t,sc \\
0  \leq s_{c,k,t}^{sc} \leq s_{a,k,t}^{sc} \quad \forall k,t,sc
\end{aligned}
\end{equation}

\begin{equation} \label{SB_3}
0  \leq L_{s,b,t}^{sc} \leq L_{b,t} \quad \forall b,t,sc \\
\end{equation}

\vspace{-1em}

\begin{equation} \label{SB_4}
-r_{d,i,t}  \leq r_{a,i,t}^{sc} \leq r_{u,i,t} \quad \forall i,t,sc \\
\end{equation}

\vspace{-1.5em}

\begin{equation} \label{SB_5}
\begin{aligned}
& \sum_{{\it i}=1}^{{\it I}}(p_{i,t}+r_{a,i,t}^{sc})+\sum_{{\it j}=1}^{{\it J}}(w_{a,j,t}^{sc}-w_{c,j,t}^{sc}) \\
& +\sum_{{\it k}=1}^{{\it K}}(s_{a,k,t}^{sc}-s_{c,k,t}^{sc})=L_{t}-\sum_{{\it b}=1}^{{\it Nb}}L_{s,b,t}^{sc} \quad \forall t,sc
\end{aligned}
\end{equation}

\vspace{-1.5em}

\begin{equation} \label{SB_6}
\begin{aligned}
\begin{aligned}
\begin{aligned}
&|\sum_{i=1}^{I}k_{l,i}(p_{i,t}+r_{a,i,t}^{sc})+\sum_{j=1}^{J}k_{l,j}(w_{a,j,t}^{sc}-w_{c,j,t}^{sc}) \\
&+\sum_{k=1}^{K}k_{l,k}(s_{a,k,t}^{sc}-s_{c,k,t}^{sc})-\sum_{{\it b}=1}^{{\it Nb}}k_{l,b}(L_{b,t}-L_{s,b,t}^{sc})|\\ & \le  Pl_l^{max} \quad \forall l,t,sc
\end{aligned}
\end{aligned}
\end{aligned}
\end{equation}
where
\begin{itemize}

	\item \eqref{SB_2} is the actual amount of REC constraint; $w_{a,j,t}^{sc}$ is actual wind power of $j$'th WPP at time {\it t} of $sc$'th scenario; $s_{a,k,t}^{sc}$ is actual solar power of $k$'th PVPP at time {\it t} of $sc$'th scenario;

	\item \eqref{SB_3} is the actual amount of LS constraint;

	\item \eqref{SB_4} is the actual amount of reserve constraint; $r_{a,i,t}^{sc}$ is actual amount of reserve of $i$'th CPP at time {\it t} of $sc$'th scenario;

	\item \eqref{SB_5} is the supply-demand balance constraint;

	\item \eqref{SB_6} is the transmission capacity constraint;


\end{itemize}


Compared with the proposed distribution-based ED model, the scenario-based ED model can not only model the cost of LS and REC caused by system reserve deficiency but also can model the cost of LS and REC caused by transmission congestion. However, the number of scenarios after reduction is limit due to the computation ability. This would reduce the representation accuracy of renewable power in the above scenario-based RTED model, which would be discussed in Section VI.

When the reserve deficiency or transmission congestion occur, REC and LS have to be employed for the balance of system power. Optimal REC and LS strategies can be obtained by solving the static optimization problem (ob. \eqref{SB_1}, s.t. {\eqref{SB_2}-\eqref{SB_6})} by the deterministic value of CPPs scheduled power, actual reserve, actual power of WPPs and PVPPs (the only scenario in this optimization problem).





\section{Case Study}

The IEEE 118-bus system is employed to validate the stochastic dynamic RTED model with multiple RPPs. There are 10 WPPs and 4 PVPPs each with a capacity of 200MW, {connecting on the 10, 24, 25, 26, 61, 65, 69, 72, 73, 87, 89, 91, 111 and 113 buses, respectively.} The data of RPPs are obtained by synchronous data in Kansas 2006 produced by NREL~\cite{Nrel}.
Their corresponding forecast power of each RPP is generating by the persistence forecast method. Gaussian Copula is used to model the conditional distribution that needed in this paper. We consider a time period with 12 intervals, modeling the 5 minute dispatch within an hour.


\subsection {Uncertainties Modeling of Multiple RPPs}

\begin{figure}[ht]
	\begin{center}
		\includegraphics[trim = 30 300 40 270, clip, width=1.0\columnwidth]{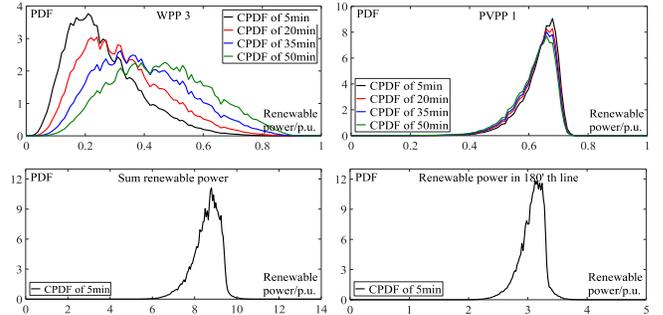}\\
		\caption{Conditional PDFs of actual outputs}\label{actual}
	\end{center}
\end{figure}

Conditional PDF of actual power of 3'th WPP, 1'th PVPP, sum renewable power and renewable power in 180'th transmission line  are as shown in Fig.~\ref{actual}. We can see that with the increase of forecast power of 3'th WPP from 5min-60min, the location of marginal distribution of actual power also moved to right. Although the conditional PDF of wind power seems to be relatively fat, the forecast error of sum renewable power tends to be thinner due to the independence of WPPs and PVPPs.

\begin{figure}[ht]
	\begin{center}
		\includegraphics[trim = 30 350 60 370, clip, width=1.0\columnwidth]{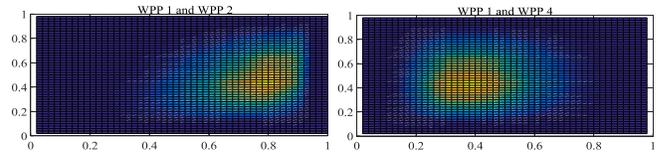}\\
		\caption{Conditional joint PDFs of actual output of two RPPs}\label{jointpdf}
	\end{center}
\end{figure}

To further analyze the renewable power independence, conditional joint PDFs of two RPPs of {\it t}=1 are shown in Fig.~\ref{jointpdf} for its high sensitivity. We can see that the power of 1'th WPP and 3'th WPP show the feature of positive correlation since they are geographically close. In contrast, the 1'th WPP and 4'th WPP are further apart and has smaller correlations.

\vspace{-1em}

\subsection {Renewable Power Scenarios}
To generate renewable power scenarios that captures the spatial correlation, renewable power scenarios should follow the joint distribution in \eqref{cjdistribution}. However, it is usually hard to use \eqref{cjdistribution} directly. For instance, a $100^{14}$ size matrix would be needed to store the joint distribution with 0.01p.u. resolution in this case. In contrast, thanks to the Gibbs theory, renewable power scenarios can be generated by the proposed method with only 100 size matrix to store the conditional distribution in \eqref{ccdistribution} with same resolution.

\begin{figure}[ht]
	\begin{center}
		\includegraphics[trim = 135 340 130 360, clip, width=1.0\columnwidth]{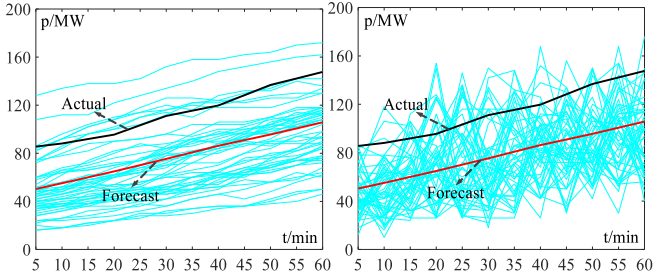}\\
		\caption{The left picture shows the scenarios generated when time correlations are considered (our method) v.s. scenarios that do not consider correlations in time (on the right, standard method). By considering temporal correlations, much more realistic scenarios can be generated.}\label{scenario1}
	\end{center}
\end{figure}

To show the effect of variability (temporal correlation), we generate renewable power scenarios by our proposed method and the method in \cite{copula_Zhang}, respectively. $N_{sc}$ is 5000, $N_{b}$ is 1000 in this case. Fig.~\ref{scenario1} shows the former 50 scenarios in $N_{sc}$ of 3'rd WPP based on our proposed method~(the left figure) and the method in \cite{copula_Zhang}~(the right figure). The red line and black line in Fig.~\ref{scenario1} is the forecast power and actual power, respectively. We can see that the above two scenarios set have the same distribution in each time interval while the renewable power scenarios of our method are much more similar to the actual renewable power. The economic comparison is discussed in Section~\ref{sec:econ}.


\vspace{-1em}

\subsection {Economy Comparison of Different RTED Methods} \label{sec:econ}

The generated renewable power scenarios are reduced to 10 scenarios and incorporated in proposed scenario-based ED. To compare the economy of the proposed RTED model, scheduled power of CPPs obtained by different RTED model are tested with other generated 10000 scenarios that consider the variability. {The cost coefficients of upward and downward reserve are all 10\$/MW.} The penalty coefficients of LS and REC is 1000\$/MW and 80\$/MW, respectively. The following five RTED models are compared in this paper.

{\it Case1}: The proposed distribution-based RTED model. {\it Case2}: The proposed distribution-based RTED model while the transmission capacity constraint use the forecast renewable power in \cite{VersatileMixture}. {\it Case3}: The proposed scenario-based RTED model. {\it Case4}: The proposed scenario-based RTED model while does not consider the variability of renewable power as in \cite{copula_Zhang}. {\it Case5}: Scenario-based RTED model that uses the marginal distribution of each RPP by \cite{sce_generation_Ma}, i.e. does not consider the spatial correlation of RPPs. The average costs of the above five cases are shown in Table.~\ref{cost}.

\begin{table}[h]
	\caption{{Total cost of different RTED models}}
	\label{cost}
	\begin{center}
		\begin{tabular}{|p{2.6cm}<{\centering}|c|c|c|c|c|}
			\hline
			{Cost/\$} & Fuel & Reserve & LS & REC & Total\\ \hline
			Proposed distribution- based model & \rule{0pt}{0.3cm} 36584 &4590&308&830& 42312\\\hline
			Model in \cite{VersatileMixture} &35877  & 4589& 702& 5705& 46873 \\\hline
			Proposed scenario- based model & \rule{0pt}{0.3cm} 35010 & 2986& 4514& 5298& 47808\\\hline
			Model in \cite{copula_Zhang} & 35032 & 2667& 7018& 5361& 50078\\ \hline
			Model in \cite{sce_generation_Ma} &35002& 2710 & 7741& 5425& 50878  \\ \hline

		\end{tabular}
	\end{center}
\end{table}

\begin{figure}[ht]
	\begin{center}
		\includegraphics[trim = 35 230 35 190, clip, width=1.0\columnwidth]{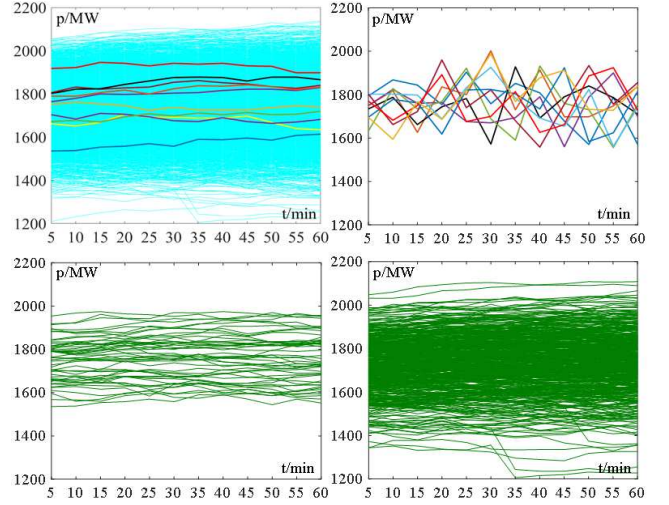}\\
		\caption{{Original scenarios and scenarios after reduction.}}\label{reducedscenarios}
	\end{center}
\end{figure}

{As shown in Table.~\ref{cost}, compared with the proposed distribution-based model, model in \cite{VersatileMixture} has larger LS and REC penalty for the reason that it has not considered the transmission congestion caused by renewable power uncertainties. Compared with the proposed scenario-based model, model in \cite{copula_Zhang} has larger cost for the reason that it could not capture the renewable power variability. Compared with the proposed scenario-based model, model in \cite{sce_generation_Ma} has larger cost for the reason that it has not considered the correlations between different RPPs. Overall, the scenario-based RTED method has much larger LS and REC penalty for the reason that it underestimate the uncertainties of renewable power compared with distribution-based RTED method, as shown in Fig.~\ref{reducedscenarios}. The slim blue lines in the upper left figure are the original renewable scenarios by the proposed scenario generation method and the figure also shows the 10 scenarios after reduction. Renewable power scenarios generated by \cite{copula_Zhang} are also reduced to 10, as shown in the upper right figure.}

{To further analyze the scenario-based RTED, the original renewable scenarios by the proposed scenario generation method are reduced to 50, 500 and 2000 scenarios and embedded in the RTED model. As shown in Table.~\ref{cost2}, with more scenarios embedded in the ED model, system cost decreases. The 50 and 500 scenarios after reduction are shown in the lower left figure and lower right figures, respectively. It can be seen that with 500 scenarios, scenario-based RTED has similar performance with the distribution-based model for the reason that renewable power uncertainties could be well represented. When 2000 scenarios embedded in the RTED model, scenario-based ED shows a better performance compared with of the distrbution-based one since it has a more flexible manner to model the risk of renewable power such as REC caused by certain transmission line congestion.}



\begin{table}[h]
	\caption{{Total cost of proposed scenario-based model based on different numbers of scenarios after reduction}}
	\label{cost2}
	\begin{center}
		\begin{tabular}{|p{2.5cm}<{\centering}|c|c|c|c|c|}
			\hline
			{Cost/\$} & Fuel & Reserve & LS & REC & Total\\ \hline
			50 scenarios& 35083 & 3124& 3464& 3954& 45625\\\hline
			500 scenarios& 36572  & 4584& 319& 854& 42329\\\hline
			2000 scenarios& 36591 &4528&285&785& 42189\\\hline
		\end{tabular}
	\end{center}
\end{table}


\vspace{-1em}

\subsection {RTED With Different Penalty Coefficients}

To show the relationships of system reserve and the potential risk of REC and LS, different penalty coefficients of LS and REC are set to change the potential risk of REC and LS. As shown in Table.~\ref{aci}, the confidence level of enough downward reserve increases if the penalty coefficients of REC increase. This means that when the potential risk of REC increase, more downward reserves are employed for the overall economy. In addition, the confidence level of enough downward reserve of different time interval in Table.~\ref{aci} varies due to the different uncertainties. Compared with other classical stochastic ED methods that use a predefined confidence level, the proposed can seek for the optimal confidence level.

\begin{table}[h]
	\caption{Confidence level of enough downward reserve under different REC penalty coefficients}
	\label{aci}
	\begin{center}
		\begin{tabular}{|c|c|c|c|c|c|}
			\hline
			{} & \multicolumn{5}{c|}{Penalty coefficients of REC (\$/MW$\cdot$h)}\\
			\hline
			{{\it time}/min} & {40} & {60} & {80} & {120} & {200}\\
			\hline
			05 & 74.84\% & 80.05\% & 81.98\% & 86.74\% & 93.05\% \\
			10 & 75.04\% & 80.26\% & 82.15\% & 87.08\% & 92.62\% \\
			15 & 74.73\% & 79.93\% & 81.82\% & 87.35\% & 92.16\% \\
			20 & 74.98\% & 81.29\% & 83.25\% & 86.88\% & 92.61\% \\
			25 & 74.99\% & 79.91\% & 82.64\% & 86.89\% & 93.11\% \\
			30 & 75.32\% & 81.92\% & 84.68\% & 88.88\% & 91.74\% \\
			35 & 74.98\% & 81.37\% & 84.12\% & 87.31\% & 92.06\% \\
			40 & 75.03\% & 81.80\% & 83.99\% & 88.67\% & 93.18\% \\
			45 & 75.14\% & 81.89\% & 84.08\% & 87.93\% & 93.23\% \\
			50 & 74.97\% & 81.75\% & 83.95\% & 88.63\% & 93.72\% \\
			55 & 74.71\% & 81.17\% & 83.98\% & 87.71\% & 92.20\% \\
			60 & 74.98\% & 80.86\% & 83.70\% & 88.69\% & 93.10\% \\
			\hline
		\end{tabular}
	\end{center}
\end{table}

\section{Conclusion}

This paper considers the uncertainties and correlations of multiple RPPs in real-time economic dispatch problems. We propose two methods, distribution-based and scenario-based dispatch models that take into account of system reserve and transmission congestion. We propose a scenario generation method that greatly reduces the required computational complexity and can accurately represent renewable power uncertainties, spatial correlation and variability. Results show that although the scenario-based RTED method has a better consideration in the effect of uncertainties and correlations on the system in theory, the discrete feature of scenarios after reduction greatly reduces the effect. Compared with other RTED models, the proposed methods show better economy by capturing renewable power uncertainties, spatial correlation and variability.
%


%










\bibliographystyle{IEEEtran}
\bibliography{Bibliography}
\end{document}